\documentclass[12pt]{article}
\usepackage{amsfonts,amssymb,latexsym,amsmath}
\newtheorem{theorem}{Theorem}[section]

\newtheorem{proposition}[theorem]{Proposition}

\newenvironment{proof}
{\par\addvspace{0.3cm}\noindent{\rm Proof. }}
{\nopagebreak\mbox{}\hfill $\Box$\par\addvspace{0.25cm}}

\newcommand{\be}{\begin{equation}}
\newcommand{\ee}{\end{equation}}

\newcommand{\bq}{\begin{eqnarray}}
\newcommand{\eq}{\end{eqnarray}}
\newcommand{\nn}{\nonumber}
\newcommand{\ba}{\begin{array}}
\newcommand{\ea}{\end{array}}
\newcommand{\comment}[1]{} %%%%
\newcommand{\tr}{\mathrm{trace}}

\newcommand{\iv}{^{-1}}
\newcommand{\iy}{\infty}

\newcommand{\Z}{{\mathbb Z}}
\newcommand{\T}{{\mathbb T}}

\newcommand{\cS}{\mathcal{S}}

\newcommand{\cC}{\mathcal{C}}

\renewcommand{\rho}{{\varrho}}
\renewcommand{\kappa}{\varkappa}

\newcommand{\cL}{{\cal L}}

\newcommand{\cT}{\mathcal{T}}

\begin{document}

\date{}
\title{Determinant computations for some classes of Toeplitz-Hankel matrices}
\author{Estelle L. Basor\thanks{ebasor@calpoly.edu.
          Supported in part by NSF Grant DMS-0500892}\\
                         Department of Mathematics\\
               California Polytechnic State University\\
               San Luis Obispo, CA 93407, USA
        \and
        Torsten Ehrhardt\thanks{ehrhardt@math.ucsc.edu}\\
              Department of Mathematics \\
          University of California \\
               Santa Cruz, CA-95064, USA}
\maketitle

\begin{abstract} The purpose of this paper is to compute the asymptotics of determinants of finite sections of operators that are trace class perturbations of Toeplitz operators. For example, we consider  the asymptotics in the case where the matrices are of the form $ (a_{i-j} \pm a_{i+j+1-k})_{i,j=0\dots N-1} $ with $k$ is fixed. We will show that this example as well as some general classes of operators have expansions that are similar to those that appear in the Strong Szeg\"{o} Limit Theorem. We also obtain exact identitities for some of the determinants that are analogous to the one derived independently by Geronimo and Case and by Borodin and Okounkov for finite Toeplitz matrices. These problems were motivated by considering certain statistical quantities that appear in random matrix theory.
\end{abstract}
\section{Introduction}
There is a fundamental connection between determinants of certain
matrices and random matrix ensembles. For example, one can
consider the Haar measure as a probability measure on the set of
$N \times N$ unitary matrices. This set of matrices along with
this measure is usually referred to as the Circular Unitary
Ensemble (CUE). From this measure one can then compute the density
for the distribution of the eigenvalues $e^{i\theta_{1}}, \dots,
e^{i\theta_{N}}$ of the matrices. The density turns out to be a
constant times
\[ \prod_{j < k}|e^{i\theta_{j}} -e^{i\theta_{k}}|^{2} .\]
A linear statistic for this ensemble is a random variable of the
form
\[ X_f=\sum_{j=1}^{N} f(e^{i\theta_{j}}),\]
and it is this quantity which is connected to a Toeplitz
determinant. More precisely,
\[
\frac{1}{(2\pi)^{N}N!}
\int_{-\pi}^{\pi} \dots \int_{-\pi}^{\pi} \prod_{j=1}^{N} e^{i
\lambda f(e^{i\theta_{j}})} \prod_{j < k}|e^{i\theta_{j}}
-e^{i\theta_{k}}|^{2} d\theta_{1} \dots d\theta_{N}
\]
is identically equal to
\[
\det \left(\frac{1}{2\pi}\int_{-\pi}^{\pi }e^{i\lambda
f(\theta)}e^{-i(j-k)\theta}d\theta \right)_{j,k=0, \dots ,N-1}.
\]
In probability terms this means that the inverse Fourier transform of the density is a Toeplitz determinant. In the opposite sense, the Toeplitz determinant can be thought of as an average or expectation with respect to CUE. For a proof of this and for more general facts about random matrices, we refer the reader to \cite{M, TW}.

Asymptotics of the determinant gives us information about the linear statistic. This is especially useful when the function $f$ is smooth enough, because we may appeal to the Strong Szeg\"{o} Limit Theorem to tell us asymptotically the behavior of the density function. For general Toeplitz determinants we consider
\[ \det \left(a_{j-k} \right) _{j,k  = 0, \dots ,N-1}\]
where $a_k$ denotes the $k$th Fourier coefficients of some
function $a\in L^1(\T)$ defined on the unit circle. Under
appropriate conditions the Szeg\"o Limit Theorem 
(see, e.g., \cite{BS, W2}) states that 
\bq\label{f.0}  
\det \left(a_{j-k} \right) _{j,k  = 0, \dots ,N-1} &\sim&  G^{N}[a]E[a]
\eq
as $N\to \iy$, where \bq\label{f.G1}
 G[a]&=& \exp\left( \frac{1}{2\pi}\int_{-\pi}^{\pi}\log
a(e^{i\theta}) \, d\theta \right),\\
 E[a] &=& \exp \left(\sum_{k=1}^{\infty}ks_{k}s_{-k} \right)
\eq with $s_{k}$ denoting the $k$th Fourier coefficient of $\log a.$ The
reader can check that in the case of linear statistics, where the
function $a$ is of the form $e^{i\lambda f},$ this implies that
the probability distributions for linear statistics are asymptotically
Gaussian ($N\to\iy$).

It is also know that different types  of random matrix ensembles lead to different classes of
determinants. If one considers, for instance,  averages for $O^{+}(2N)$,
the set of orthogonal matrices with determinant equal to one (or.
equivalently, the density function for certain linear statistics),
then the corresponding determinant is that of a finite Toeplitz plus
Hankel matrix. More specifically, it is of the form
\[ \det \left(a_{j-k} + a_{j+k} \right) _{j,k  = 0, \dots ,N-1}\]
where the function $a$ is assumed to be even. Because of this reason we are
interested in the determinants of a sum of a finite Toeplitz plus
a ``certain type'' of Hankel matrix. For other ensembles other
determinants arise. Two other cases of interest are
\[ \det \left(a_{j-k} + a_{j+k+1} \right) _{j,k  = 0, \dots ,N-1}\]
and
\[ \det \left(a_{j-k} - a_{j+k+2} \right) _{j,k  = 0, \dots ,N-1}.\]
We refer the interested reader to \cite{BR, FF} for derivations of the
averages in the above cases and applications of the results in random matrix theory.

Our goal is to extend as much as possible the Szeg\"{o} Limit
Theorem to these various types of determinants for both smooth and
singular symbols. In this paper we address the case of smooth
symbols. An outline of the paper is as follows. In the next
section we present some Banach algebra preliminaries and compute
some operator determinants. Then we give some explicit examples of
our general theory which correspond to the ones discussed earlier.
We then return to the general setting and derive a
Borodin-Okounkov-Geronimo-Case identity for the various classes of
operators and then compute the analogue of the Szeg\"{o} Limit
Theorem. This will allow us to calculate the asymptotics of
determinants of the form
\[ (a_{i-j} \pm a_{i+j+1-k})_{i,j=0\dots N-1} \]
where $k\in\Z$ is fixed. Finally, in the last section, we present some additional results about our classes
of operators.

%%%%%%%%%

%%%%%%%%%

\section{Computation of operator determinants}
\label{s2}

We denote by $\ell^{2}$ the space of all complex-valued square-summable sequences
$\{x_{n}\}_{n=0}^{\infty}$. The set $\cL(\ell^{2})$ is the set of all 
bounded linear operators on  $\ell^{2}$. By $\cC_1(\ell^2)$ we denote the
class of trace class operators on $\ell^2$. We refer to \cite{GK} for more
information about this class and the related notions of operator
traces and determinants.

Let $a \in L^{\infty}(\T)$ be a measurable and essentially bounded
function on the unit circle. The Toeplitz operator $T(a)$ and
Hankel operator $H(a)$ with symbol $a$ are the bounded linear
operator defined on $\ell^{2}$ with matrix representations
$$
T(a) = (a_{j-k}),\,\,\,\,\, 0\leq j,k < \infty,
$$
and
$$
H(a) = (a_{j+k+1}),\,\,\,\,\, 0\leq j,k < \infty.
$$
%We can also equivalently define $T(a)$  and $H(a)$ by the formulas
%\[ T(a)f = P(af),\qquad H(a)f=P(a (Jf)) \]
%where $(Jf)(t)=t\iv f(t\iv)$ and $P$ is the Riesz projection on $L^2(\T)$ defined by
%\[ P : \sum_{k=-\infty}^{\infty}a_{k}e^{ik\theta} \rightarrow \sum_{k = 0}^{\infty}a_{k}e^{ik\theta} .\]
%This last definition uses the standard identification of a function $f\in L^{2}(\T)$ with its sequence of Fourier coefficients. 
It is well-known that Toeplitz and Hankel operators satisfy the fundamental identities
\bq\label{T1} 
T(ab) &=& T(a)T(b) + H(a)H(\tilde b)
\\
\label{H1}
H(ab) &=& T(a)H(b) + H(a)T(\tilde b) .
\eq
In the last two identities $\tilde b (e^{i\theta}) = b(e^{-i\theta}) .$
It is worthwhile to point out that these identities imply that
\be\label{TH.1}
T(abc)=T(a)T(b)T(c),\qquad H(ab\tilde{c})=T(a)H(b)T(c)
\ee
for $a,b,c \in L^\iy(\T)$ if $a_n=c_{-n}=0$ for all $n>0$.

The Riesz projection acting  on $L^p(\T)$ ($1<p<\iy$) is defined by
\[ P : \sum_{k=-\infty}^{\infty}a_{k}e^{ik\theta} \rightarrow \sum_{k = 0}^{\infty}a_{k}e^{ik\theta} .\]

Let $\cS$ stand for a unital Banach algebra of functions on the
unit circle which is continuously embedded into $L^\iy(\T)$ and which
has the following properties:
\begin{itemize}
\item[(i)]
the Riesz projection $P: \cS\to\cS$ is well-defined and bounded on $\cS$,
\item[(ii)]
the symmetric flip $a\in \cS\mapsto \tilde{a}\in\cS$ is well
defined and bounded on $\cS$.
\end{itemize}
Then we can define
\bq
\cS_\pm &=& \Big\{ a\in \cS\;:\; a_n=0 \mbox{ for  all } (\pm n)< 0\Big\},\nn\\
\cS_0 &=& \Big\{ a\in \cS\;:\; a=\tilde{a} \Big\}.\nn
\eq
 Moreover, we can make the following basic observations.
Each $a\in\cS$ can be decomposed into $a=a_++a_-$ with $a_\pm \in
\cS_\pm$. The decomposition can be made unique by requiring  that
$[a_-]_0=0$. Then the mappings $a\mapsto a_\pm $ are linear and
bounded.

Furthermore, each $a\in \cS$ can be decomposed into $a=a_0+a_-$
with $a_0\in\cS_0$ and $a_-\in \cS_-$. Indeed, this decomposition
can be derived from the previous one by writing
$$
a=a_++a_-=(a_++\tilde{a}_+)+(a_--\tilde{a}_+)
$$
and taking $a_0=a_++\tilde{a}_+$ and $a_--\tilde{a}_+$ as the new
$a_-$. Again, we can make the decomposition  unique by requiring
$[a_-]_0=0$. The corresponding projections are bounded.

Now let $\cS$ be such a Banach algebra, and assume that $M: a\in
\cS \mapsto M(a)\in \cL(\ell^2)$ is a linear and continuous map
such that the following conditions are fulfilled:
\begin{enumerate}
\item[(a)]
If $a\in \cS$, then $M(a)-T(a)\in \cC_1(\ell^2)$ and
$$
\|M(a)-T(a)\|_{\cC_1(\ell^2)} \le C\, \|a\|_{\cS}.
$$
\item[(b)]
If $a\in \cS_-$, $b\in \cS$, $c\in \cS_0$, then
$$
M(abc)=T(a)M(b) M(c).
$$
\item[(c)]
$M(1)=I$.
\end{enumerate}
We will refer to $a$ as the symbol of $M(a)$. The pair $[M,\cS]$ will be called a {\em compatible pair}.

Let us remark that, assuming (c), condition (b) is equivalent to the conditions that
\be\label{f.two}
M(ab)=T(a)M(b), \quad M(bc)=M(b)M(c)
\ee
whenever $a\in \cS_-$, $b\in \cS$, $c\in \cS_0$.

\begin{proposition}\label{p2.1}
Let $[M,\cS]$ be a compatible pair. Then
\begin{enumerate}
\item [(i)] $H(a)H(b) \in \cC_{1}(\ell^{2})$  for each $a,b\in \cS$, and there is a constant $C$ such that
$$
\|H(a)H(b)\|_{\cC_1(\ell^2)}\le  C\|a\|_{\cS}\|b\|_{\cS}
\quad \mbox{for each } a,b\in\cS,
$$
\item [(ii)]
if $a$ is invertible in $\cS$, then
$T(a\iv)M(a) -I$ and $M(a)T(a\iv) - I$ are both in $\cC_{1}(\ell^{2}).$
\end{enumerate}
\end{proposition}
\begin{proof}
(i): We first assume that $b = \tilde b.$
By assumption (a) each of the operators
\[ (M(a) - T(a) )M(b), \,\,\,\,
T(a)(M(b) - T(b) ), \,\,\,\,T(ab)-M(ab) \] is in
$\cC_{1}(\ell^{2}).$ If we add these three operators together
and use that $M(ab)=M(a)M(b)$ by proporty (b), it follows that
\[ T(ab) - T(a)T(b)=H(a)H(b) \]
is trace class. With a more careful inspection we can furthermore derive the estimate
$$
\|H(a)H(b)\|_{\cC_1(\ell^2)}\le C \|a\|_{\cS}\|b\|_{\cS}.
$$
In general we write $b=b_0+b_-$ with $b_0\in\cS_0$ and $b_-\in\cS_-$.
Then $H(b)=H(b_0)$ and the result follows. The estimate also holds because the map $b\mapsto b_0$ is bounded.

(ii): We note that if $a$ is
invertible, then
\[ T(a^{-1}) ( M(a) - T(a) ) \] is trace class and by the first  part \[T(a^{-1}) T(a) - I = -H(a^{-1})H(\tilde a)\] is also trace class. The proof for $M(a)T(a\iv) - I$ is similar.
\end{proof}

Statement (i) of the previous proposition implies that if
$[M,\cS]$ is a compatible pair, then $\cS$ is a suitable Banach
algebra in the sense of \cite{E4}. It has been shown there that in such a setting the 
the Szeg\"o Limit Theorem, i.e., the asymptotics (\ref{f.0}),
holds for each $a\in \mathcal{G}_1 \cS$. Here $\mathcal{G}_1 \cS$
stands for the group of exponentials of functions in $\cS$.

In the following proposition we employ the notion of an analytic
Banach algebra valued function. We refer to \cite{Ru} for details.

\begin{proposition}
Let $[M,\cS]$ be a compatible pair. Then for $a\in \cS$, the functions
$$
F(\lambda)=T(e^{-\lambda a}) M(e^{\lambda a})-I
$$ and
$$
F(\lambda)=T\iv(e^{\lambda a}) M(e^{\lambda a})-I
$$
are $\cC_1(\ell^2)$-valued analytic.
\end{proposition}
\begin{proof}
The first function is obviously an analytic $\cL(\ell^2)$-valued function. It is trace class valued because of
(ii) of the previous proposition.

In order to consider the second function, 
we decompose $a=a_++a_-$ with $a_\pm \in\cS_\pm $.
From this we derive that
\[T(e^{\lambda a}) = T(e^{\lambda a_{-}})T(e^{\lambda a_{+}}).
\]
Hence the inverse exists and is given by
\be\label{f.1}
T(e^{\lambda a})\iv = T(e^{-\lambda a_{+}})T(e^{-\lambda a_{-}}).
\ee
This shows that also the second function is well defined and $\cL(\ell^2)$-valued analytic. By assumption (a),
it is trace class valued.
\end{proof}

We now compute some operator determinants in the next two propositions. They will appear later as 
constants in our asymptotic relations. We use the following well known facts. If
$F(\lambda)$ is an analytic function of the form identity plus trace class, then determinant $\det F(\lambda)$ is well defined and (complex-valued) analytic. Moreover, 
$$
(\log \det F(\lambda))'=\frac{(\det F(\lambda))'}{\det F(\lambda)}=\tr \,F'(\lambda) F\iv(\lambda)=
\tr\,  F\iv(\lambda) F'(\lambda).
$$
The proof of the following propositions is similar to the proof of, for instance, \cite[Thm.~2.5]{BE1} and \cite[Thm.~7.4]{E4}.

\begin{proposition}\label{p2.3}
Let $[M,\cS]$ be a compatible pair. Then for $a\in \cS_0$,
\be\label{1}
\det T (e^{-a})M(e^a)= \exp\Big( \tr(M(a)-T(a)) + \frac{1}{2}\tr\, H(a)^2\Big).
\ee
\end{proposition}
\begin{proof}
Define the entire function
$$
f(\lambda):=\det T(e^{-\lambda a} )M(e^{\lambda a}).
$$
Now consider the logarithmic derivative of $f(\lambda)$,
\bq
\frac{f'(\lambda)}{f(\lambda)} & =&
\tr \Big ( M(e^{-\lambda a}) T^{-1}(e^{-\lambda a} )\Big)\Big(T(e^{-\lambda a}) M(ae^{\lambda a}) - T(a e^{-\lambda a})M(e^{\lambda a})\Big)\nn\\
&=&
\tr \Big(  M(a) -T\iv(e^{-\lambda a}) T(a e^{-\lambda a})
\Big).\nn
 \eq
Differentiating again yields
\bq
\left(\frac{f'(\lambda)}{f(\lambda)}\right)' &=&
\tr \Big(-T\iv(e^{-\lambda a}) T(a e^{-\lambda a })T\iv(e^{-\lambda a})T(ae^{-\lambda a}) +
T^{-1}(e^{-\lambda a})T(a ^2e^{-\lambda a})
\Big)\nn\\
&=& \tr\Big(-T(a )T(a)+T(a^2)\Big)=\tr\, H(a)^2.\nn
\eq
The last equality holds by writing $a = a_- + a_{+}$ with $a_\pm\in \cS_\pm$ and
considering the inverse of $T(e^{-\lambda a})$ as in (\ref{f.1}).
Integration and fixing the constants by putting $\lambda=0$ yields
\bq
f(\lambda)=\exp\Big(\lambda\, \tr(M(a)-T(a))+\frac{\lambda^2}{2}\tr\, H(a)^2\Big).\nn
\eq
This finishes the proof.
\end{proof}

\begin{proposition}\label{p2.4}
Let $[M,\cS]$ be a compatible pair. Then for $a\in \cS$,
\be%\label{1}
\det T\iv (e^a)M(e^a)= \exp\Big( \tr(M(a)-T(a))-\frac{1}{2}\tr\, H(a)^2\Big).
\ee
\end{proposition}
\begin{proof}
We can decompose $a=a_0+a_-$ with $a_-\in\cS_-$, $a_0\in\cS_0$.
It is easy to see that
$$
f(\lambda):=\det T\iv (e^{\lambda a} )M(e^{\lambda a})=\det T\iv (e^{\lambda a_0})M(e^{\lambda a_0}).
$$
Now consider the logarithmic derivative of $f(\lambda)$,
\bq
\frac{f'(\lambda)}{f(\lambda)}
&=&
\tr \Big( M(e^{-\lambda a_0}) T(e^{\lambda a_0})\Big)
\Big(T\iv (e^{\lambda a_0}) M(a_0e^{\lambda a_0})-T\iv(e^{\lambda a_0}) T(a_0e^{\lambda a_0})
T\iv(e^{\lambda a_0})M(e^{\lambda a_0})
\Big)\nn\\
% &=&
%\tr \Big( T\iv (e^{\lambda a_0}) M(a_0) T(e^{\lambda a_0})-T\iv(e^{\lambda a_0}) T(a_0e^{\lambda a_0})
%\Big)\nn\\
&=&
\tr \Big( M(a_0) -T(a_0e^{\lambda a_0})T\iv(e^{\lambda a_0})
\Big).\nn
\eq
Differentiating again yields
\bq
\left(\frac{f'(\lambda)}{f(\lambda)}\right)' &=&
\tr \Big(T(a_0e^{\lambda a_0})T\iv(e^{\lambda a_0}) T(a_0e^{\lambda a_0})T\iv(e^{\lambda a_0}) -
T(a_0^2e^{\lambda a_0})T\iv(e^{\lambda a_0})
\Big)\nn\\
&=& \tr\Big(T(a_0)T(a_0)-T(a_0^2)\Big)=-\tr\, H(a_0)^2.\nn
\eq
Integration and fixing the constants by putting $\lambda=0$ yields
\bq
f(\lambda)=\exp\Big(\lambda\, \tr(M(a_0)-T(a_0))-\frac{\lambda^2}{2}\tr\, H(a_0)^2\Big).\nn
\eq
This implies the desired assertion by noting that $H(a_-)=0$ and $M(a_-)=T(a_-)$ by parts
(b) and (c) of the assumptions.
\end{proof}

%%%%%%%%

%%%%%%%%

\section{Concrete realizations}
\label{s3}

While the above formulas are nice, it remains to show that there are some interesting classes of operators that satisfy the Banach algebra conditions of the previous section as well as the algebraic conditions, that is, we need to show that there are some compatible pairs. We would also like of course to have operators that correspond to the random matrix examples that were stated in the introduction. The purpose of this section is to introduce these examples, or concrete realizations. We need to specify the Banach algebra $\cS$ and to identify the operators $M(a).$

For our compatible pairs, it is convenient to take as Banach algebra the Besov class $B_{1}^{1}.$ This is the algebra of all functions $a$ defined on the unit circle for which
\[
\|a\|_{B_{1}^{1}} := \int_{-\pi}^{\pi}\frac{1}{y^{2}}\int _{-\pi}^{\pi}|a(e^{ix+iy}) + a(e^{ix-iy}) -2a(e^{ix})\,|dxdy < \iy .
\]
A function $a$ is in this class if and only if the Hankel operators $H(a)$ and $H(\tilde{a})$ are both trace class. Moreover the Riesz projection is bounded on this class and an equivalent norm is given by
\[ |a_{0}| + \|H(a)\|_{\cC_1} + \|H(\tilde{a})\|_{\cC_1} .\]
A proof of these facts can be found in \cite{Pe1, Pe2}.

Introduce the projections
$$
P_1=\mathrm{diag}(1,0,0,\dots),\qquad Q_1=I-P_1
$$
acting on $\ell^2$.
 Next we define four realizations of $M(a)$. Given our conditions on the Banach algebra $B_{1}^{1}$, we need  only check that the required properties for $M$ hold.

\begin{proposition}\label{p3.1}
 The following realizations for the operator $M$ with symbols in the Besov class $B_1^1$ define compatible pairs $[M, B_1^1]$:
\begin{enumerate}
\item[(I)]
$M(a) =T(a)+H(a)$
\item[(II)]
$M(a) =T(a)-H(a)$
\item[(III)]
$M(a)=T(a)-H(t\iv a)$
\item[(IV)]
$M(a)=T(a)+H( t a)Q_1$
\end{enumerate}
\end{proposition}
\begin{proof}
It is easily seen (by aking into account of the remark made in connection with (\ref{f.two})),
that the only not immediately obvious issue is the multiplicative property
$$M(ab)=M(a)M(b)$$ under the condition $b=\tilde{b}$.
In order to verify the cases (I) and (II) use (\ref{T1}) and (\ref{H1}) to obtain
\bq
T(ab)\pm H(ab) & =& 
T(a)T(b)+H(a)H(b)\pm T(a)H(b)\pm H(a)T(a)\nn\\
&=&
\Big(T(a) \pm H(a)\Big) \Big(T(b)\pm H(b)\Big)
\nn
\eq
as desired. In case (III), use in addition (\ref{TH.1}) to obtain 
\bq
M(a)M(b) &=& \Big(T(a)-H(t\iv a)\Big)\Big(T(b)-H(t\iv b)\Big)\nn\\
&=& T(a)T(b)+H(t\iv a) H(t\iv b)- H(t\iv a)T(b)-T(a)H(t\iv b)\nn\\
&=& T(ab)-H(a)P_1H(b)-H(t\iv a b) +T(t\iv a)H(b)-T(a)T(t\iv)H(b)\nn\\
&=& T(ab)-H(t\iv ab)-H(a)P_1 H(b) +H(a)H(t)H(b)\nn\\
&=& T(ab)-H(t\iv ab).\nn
\eq
Here $P_1=H(t)=H(t)^2=I-Q_1$ and $Q_1=T(t)T(t\iv)$.
In case (IV) we have 
\bq
M(a)M(b) &=& \Big(T(a)+H(t a)Q_1\Big)\Big(T(b)+H(t b)Q_1\Big)\nn\\[.5ex]
&=&
T(a)T(b) +H(ta) Q_1 H(t b)Q_1+
H(ta) Q_1 T(b) + T(a) H(t b) Q_1
\nn\\[.5ex]
&=&
T(a) T(b) +H(a) H(b)Q_1+
H(a)T(t\iv b)+T(a)H(tb)Q_1
\nn\\[.5ex] 
&=&
T(ab)-H(a)H(b)P_1+H(a)T(t\iv b)+H(tab)Q_1-H(a)T(t\iv b) Q_1
\nn\\[.5ex]
&=&
T(ab)+H(a)T(b) T(t\iv)+H(tab)Q_1-H(a)T(t\iv b) Q_1
\nn\\[.5ex]
&=&
T(ab)+H(ab)Q_1.\nn
\eq
This settles the proof.
\end{proof}

Let us remark that the operators (I)-(III) are precisely the infinite matrix  versions of the finite Toeplitz plus Hankel matrices  mentioned in the introduction. It is also easily seen that if we multiply 
the operator (IV) from the right with $\mathrm{diag}(2,1,1,\dots)$, then we obtain
$T(a)+H(ta)$. Finally notice the simple fact that the operators (I) and (II) are related with one another by multiplying from the left and right with $\mathrm{diag}(1,-1,1,-1,\dots)$ and replacing the symbol
$a(t)$ by $a(-t)$.

\begin{proposition}\label{p3.2}
Let $a\in B^1_1$ and denote
\bq\label{def.F}
F[a]& =& \det T\iv(a) M(a)
\eq
where we assume that there exists a logarithm $\log a\in B_1^1$.
Then in the above cases (I)--(IV) the corresponding constants evaluate as follows:
\bq
F_I[a] &=& \exp\Big(\sum_{n=0}^\iy [\log a]_{2n+1}-\frac{1}{2}\sum_{n=1}^\iy n  [\log a]_n^2\Big)
\nn\\
F_{II[}a] &=& \exp\Big(-\sum_{n=0}^\iy [\log a]_{2n+1}-\frac{1}{2}\sum_{n=1}^\iy n [\log a]_n^2\Big)
\nn\\
F_{III}[a] &=& \exp\Big(-\sum_{n=1}^\iy [\log a]_{2n}-\frac{1}{2}\sum_{n=1}^\iy n  [\log a]_n^2\Big)
\nn\\
F_{IV}[a] &=& \exp\Big(\sum_{n=1}^\iy [\log a]_{2n}-\frac{1}{2}\sum_{n=1}^\iy n  [\log a]_n^2\Big)
\nn
\eq
\end{proposition}
\begin{proof}
We only need to note that $\tr\, H(\log a)^2$ is $\sum_{n=1}^\iy n  [\log a]_n^2$ and check that for example
$\tr H(\log a) = \sum_{n=0}^\iy [\log a]_{2n+1}.$
\end{proof}
The proof of the following proposition is almost the same as above.

\begin{proposition}\label{p3.3}
Let $a\in B^1_1$ and denote
\bq\label{def.Fhat}
\hat{F}[a] &=& \det T(a\iv) M(a)
\eq
where we assume that there exists a logarithm $\log a\in B^1_1$ and $a=\tilde{a}$.
Then in the above cases (I)--(IV) the corresponding constants evaluate as follows:
\bq
\hat{F}_I[a] &=& \exp\Big(\sum_{n=0}^\iy [\log a]_{2n+1}+\frac{1}{2}\sum_{n=1}^\iy n  [\log a]_n^2\Big)
\nn\\
\hat{F}_{II[}a] &=& \exp\Big(-\sum_{n=0}^\iy [\log a]_{2n+1}+\frac{1}{2}\sum_{n=1}^\iy n [\log a]_n^2\Big)
\nn\\
\hat{F}_{III}[a] &=& \exp\Big(-\sum_{n=1}^\iy [\log a]_{2n}+\frac{1}{2}\sum_{n=1}^\iy n  [\log a]_n^2\Big)
\nn\\
\hat{F}_{IV}[a] &=& \exp\Big(\sum_{n=1}^\iy [\log a]_{2n}+\frac{1}{2}\sum_{n=1}^\iy n  [\log a]_n^2\Big)
\nn
\eq
 \end{proposition}

%%%%%%%%%%%%%%%%%%%%%%%%

\section{Exact identities for some determinants}

In this section we establish some exact identities for the finite sections of the operators considered in the previous section. These are of the Borodin/Okounkov/Geronimo/Case type and with these the asymptotics of the determinants will easily follow. For the Toeplitz analogue of this theorem see \cite{BO}.
We define the projection $P_{N}$  by
\[ 
P_{N} : \{x_n\}_{n=0}^\iy\in \ell^2 \mapsto \{y_n\}_{n=0}^\iy\in \ell^2,\qquad
y_n=\left\{\ba{cc} x_n & \mbox{if } n<N\\ 0 & \mbox{if } n\ge N\ea\right.
\]
and put $Q_N=I-P_N$.
We are interested in the determinants (where the matrices or operators are always thought of as acting on the image of the projection of the appropriate space) of
\[ P_{N} M(a) P_{N} .\]
We first take the case of even $a$.
Recall the definition of the constant $G[a]$ given in (\ref{f.G1}).

\begin{proposition}
Let $[M,\cS]$ be a compatible pair, and let $b_+\in \cS_+$. Put $a=a_+\tilde{a}_+=\exp(b)$ with  $a_+=\exp (b_+)$, $b=b_++\tilde{b}_+$.
Then
\[ \det P_{N} M(a) P_{N}  = G[a]^{N} \hat{F}[a] \det (I + Q_{N} K Q_{N} ),\]
where
\[ \hat{F}[a] = \det T(a\iv) M(a)= \exp\Big( \tr(M( b)-T(b)) + \frac{1}{2}\tr\, H(b)^2\Big),\]
and $K = M(a_{+}\iv) T(a_{+}) - I.$
\end{proposition}
\begin{proof}
We can write
\begin{align}
P_{N} M(a) P_{N} & =   P_{N} M(a) P_{N}\nn \\
& =  P_{N} T( a_{+}) T( a_{+}\iv)M(a) T( {\tilde a_{+}}) T({\tilde a_{+}}\iv)P_{N}\nn\\
& = P_{N} T( a_{+}) P_{N}T( a_{+}\iv)M(a) T( {\tilde a_{+}}\iv)P_{N} T({\tilde a_{+}})P_{N}.\nn
\end{align}
The last fact follows since for Toeplitz operators
\[    P_{N} T( a_{+}) =  P_{N} T( a_{+}) P_{N},\qquad
 T({\tilde a_{+}})P_{N} = P_{N} T({\tilde a_{+}})P_{N}.\]
At this point we have that
\begin{align}
 \det P_{N} M(a) P_{N} & = \det (P_{N} T( a_{+}) P_{N}T( a_{+}\iv)M(a) T( {\tilde a_{+}}\iv)P_{N} T({\tilde a_{+}})P_{N})\nn\\
 & = \det (P_{N} T( a_{+}) P_{N} )\cdot \det (P_{N} T( {\tilde a_{+}}) P_{N}) \cdot \det P_{N}T( a_{+}\iv)M(a) T( {\tilde a_{+}}\iv)P_{N} \nn\\
  & = [a_{+}]_{0}^{N}\cdot [{\tilde a_{+}}]_{0}^{N}\cdot \det P_{N}T( a_{+}\iv)M(a) T( {\tilde a_{+}}\iv)P_{N}.
  \nn
  \end{align}
First it is not hard to check that
$ [a_{+}]_{0}\cdot [{\tilde a_{+}}]_{0} = G[a]$.
Now Jacobi's identity for invertible operators on Hilbert space which are of the form identity plus trace class operators states that for projections $P$ and $Q=I-P$ we have
\[ \det\,PAP=(\det\,A)\cdot(\det\,QA\iv Q).\label{detdet}\]
We apply this to the above with $P=P_{N }$, $Q = I - P_{N}$,
and $A=T(a_+\iv)M(a)T(\tilde{a}_+\iv)$
to find that
\[ \det P_{N}T( a_{+}\iv)M(a) T( {\tilde a_{+}}\iv)P_{N} = \det T( a_{+}\iv)M(a) T( {\tilde a_{+}}\iv)
\cdot \det Q_{N}
(T( a_{+}\iv)M(a) T( {\tilde a_{+}}\iv))\iv Q_{N}.\]
To simplfy the last two determinants we note that
\[ \det T( a_{+}\iv)M(a) T( {\tilde a_{+}}\iv) = \det T( {\tilde a_{+}}\iv)T( a_{+}\iv)M(a)  = \det T(a\iv) M(a)
=\hat{F}[a]\]
and use Proposition \ref{p2.3}. Moreover,
\[ (T( a_{+}\iv)M(a) T( {\tilde a_{+}}\iv))\iv = T( {\tilde a_{+}}) M(a\iv) T(a_{+})=M(a_+\iv)T(a_+)\]
which is also of the form $I$ plus a trace class operator.
We now put all these together,
\[  \det P_{N} M(a) P_{N}  = G[a]^{N} \hat{F}[a] \det Q_{N}  M(a_{+}\iv) T(a_{+} )Q_{N} ,\]
and make the observation that this last determinant is the same as
\[\det (P_{N}+ Q_{N}  M(a_{+}\iv) T(a_{+} )Q_{N} )
=\det (I + Q_{N}  K  Q_{N}). \]
This proves the formula.
\end{proof}

Let us remark that the operator $K$ appearing in the previous proposition becomes particularly
simple in the cases of the four concrete realizations of operators $M$ considered in the previous section.
The precise expressions are as follows:
\begin{itemize}
\item[(I)] $K=H(a_+\iv \tilde{a}_+)$
\item[(II)] $K=-H(a_+\iv \tilde{a}_+)$
\item[(III)] $K=-H(t\iv a_+\iv \tilde{a}_+)$
\item[(IV)] $K=H(t a_+\iv  \tilde{a}_+)-T(a_+\iv)H(t\tilde{a}_+)$
\end{itemize} 
As for case (IV) notice that the term $T(a_+\iv)H(t\tilde{a})$ will be annihilated by multiplication with $Q_N$ from the right ($N\ge1$).

The above proposition needs to be slightly changed for non-even functions $a$. We include the result for completeness sake, although in our applications $a$ is always even.

\begin{proposition}
Let $[M,\cS]$ be a compatible pair, and let $b_\pm\in \cS_\pm$. Put $a=a_+a_-=\exp(b)$ with  $a_\pm=\exp (b_\pm)$, $b=b_++b_- $.
Then
\[ \det P_{N} M(a) P_{N}  = G[a]^{N} E[a] F[a] \det (I + Q_{N} K Q_{N} ),\]
where
\bq
F[a] &=& \det T\iv(a) M(a) = \exp\Big( \tr(M( b)-T(b)) - \frac{1}{2}\tr\, H(b)^2\Big)
\nn\\
E[a] &=& \det T(a\iv) T(a)= \exp\Big( \tr \, H(b)H({\tilde b})\Big) ,\nn
 \eq
and $K = M(a_-a_+\iv\tilde{a}_{+}\iv) T(a_{+}\tilde{a}_+ a_-\iv) - I.$
\end{proposition}
\begin{proof}
The only real difference in the proof is we must replace $\tilde{a}_+$ by $a_-$.
The subsequent computations must be modified as follows.
Firstly,
$$
\det T(a_+\iv) M(a) T(a_-\iv)=\det T(a_-\iv)T(a_+\iv) M(a)
=\det T(a\iv) M(a).
$$
This we can write as the product
$$
\det T(a\iv) T(a)\cdot \det T(a)\iv M(a)
$$
and use Proposition \ref{p2.4} to identify the first factor. The second factor is well known 
from Toeplitz theory \cite{BS, W2}. Since in our setting the Banach algebras $\cS$
is not specified, one way to settle the issue is to use the same ideas as in Propositions
\ref{p2.3} and \ref{p2.4}. Another possibility would  be to apply a formula due to Pincus
(see, e.g., \cite{E3}).

Secondly, one computes that
\begin{align}
(T(a_+\iv)M(a) T(a_-\iv))\iv &=(T(a_+\iv)T(a_-\tilde{a}_+\iv )M(a_+\tilde{a}_+) T(a_-\iv))\iv
\nn\\
&= T(a_-) M(a_+\iv \tilde{a}_+\iv) T(a_-\iv \tilde{a}_+) T(a_+)
\nn\\
&= M(a_-a_+\iv \tilde{a}_+\iv) T(a_-\iv \tilde{a}_+ a_+).\nn
\end{align}
This is again identity plus trace class by Proposition \ref{p2.1}.
\end{proof}

The above proposition immediately establishes an asymptotic formula for the determinants since
the operators $Q_{N}$ tend to zero strongly as does its adjoint. Thus we arrive at final
results of this section.

\begin{theorem}\label{t4.3}
Let $[M,\cS]$ be a compatible pair, let $b\in \cS$ and $a=\exp(b)$.
Then
 \[ \det P_{N} M(a) P_{N}  \sim  G[a]^{N} \hat{E}[a]
 \quad \mbox{ as }N\to \iy,\]
where
\[ 
\hat{E} [a] = \exp\Big( \tr(M( b)-T(b)) - \frac{1}{2}\tr\, H(b)^2 + \tr \, H(b)H({\tilde b})\Big).\]
\end{theorem}

It is clear that the formula for $E[a]$ in the previous two theorem simplifies if
$b$ is assumed to be even and then correspond to Proposition 4.1. Also, in the case of the concrete realizations the traces of $M(a) - T(a)$ are explicit in Proposition 3.2..

%\begin{theorem} Suppose  $a$ is invertible, is even, and is in the Besov class. Then
 %\[ \det P_{N} M(a) P_{N}  \sim  G[a]^{N} E(a)\]
 %\[ E(a) = \exp\Big( \tr(M( b)-T(b)) + \frac{1}{2}\tr\, H(b)^2\Big).\]
%\end{theorem}

We remark here that in the examples of our concrete realizations and with the symbol $e^{i \lambda f}$ these theorems tell us that the distributions of the linear statistics are all Gaussian since as a function of $\lambda$ the transforms are exponentials of quadratic functions.

%{\bf P.S} We can give the expressions for $K$ in cases (I)-(VI) at  least in the even case.

We end this section with an application of the above asymptotics which yield an expansion for determinants of finite sections of operators of the form
$T(a)  \pm  H(at^{k}).$ These operators are (for general $k$)  not the ones that yield a compatible pair realization, but using Jacobi's identity
\be\label{f.J}
\det\,PAP=(\det\,A)\cdot(\det\,QA\iv Q),
\ee
we can still compute the determinants of their finite sections asymptotically. 
We prepare with a basic result still relating to compatible pairs.

\begin{proposition} 
Let $[M,\cS]$ be a compatible pair. Suppose that $a\in\cS$ such that $\log a\in\cS$. Then there exists a unique factorization of the form
\bq\label{f.fac}
a(t) &=& a_-(t) a_0(t)
\eq
such that $a_-,a_-\iv \in \cS_-$, $a_0,a_0\iv\in\cS_0$, and $[a_-]_0=0$. Moreover,
$M(a)$ is invertible and 
\bq\label{M.iv}
M(a)\iv &=& M(a_0\iv) T(a_-\iv).
\eq
\end{proposition}
\begin{proof}
We can decompose $\log a$ into $b_-+b_0$  with $b_-\in\cS_-$, $b_0\in\cS_0$, and $[b_-]_0=0$.
Then we simply put
$a_-=\exp(b_-)$ and $a_0=\exp(b_0)$ in order to obtain the factorization. 
Notice that $S_-$ and $S_0$ are unital Banach subalgebras of $\cS$. To obtain the uniqueness of the factorization write
$$
a=a_-^{(1)}a_0^{(1)}=a_-^{(2)}a_0^{(2)}
$$
whence $(a_-^{(2)})\iv a_-^{(1)}=a_0^{(2)}(a_0^{(1)})\iv$. Apply the fact that the intersection of
$\cS_-$ and $\cS_0$ are the constant functions only and use the normalization condition to 
conclude that the last products are in fact equal to one.

Using the factorization and the basic properties of $M(a)$ we have
$M(a)=T(a_-)M(a_0)$ and
$$
T(a_-)T(a_-\iv)=T(a_-\iv)T(a_-)=I,\quad M(a_0)M(a_0\iv)=M(a_0\iv)M(a_0)=I.
$$
Hence the invertibility of $M(a)$ follows.
\end{proof}

\begin{theorem}
Let $ a\in B_1^1$ such that $\log a\in B_1^1$. Assume that $a=a_-a_0$ is a factorization
of the form (\ref{f.fac}). 
\begin{enumerate}
\item
Suppose that $k$ is a negative even integer ($k=-2l$, $l\ge1$).
Then
\bq
\det P_{N}  (T(a) \pm H(at^{k}))P_{N}  &\sim&  
G[a]^{N+l} E_{1,\pm}[a] \det P_{l}(T(a_0^{-1}) \pm H(a_0^{-1}) )P_{l}\nn
\eq
as $N\to\iy$, where
\[ 
E_{1,\pm}[a] = \exp\Big(\pm \sum_{n =1}^{\infty} \log a_{2n+1}  -\frac{1}{2}\sum_{n=1}^\iy n  [\log a]_n^2 +\sum_{n=1}^\iy n [\log a]_{-n} [\log a]_n\Big).
\]
\item
Suppose that $k$ is a negative odd integer less than $-1$ ($k=-1-2l$, $l\ge 1$).
Then
\bq
\det P_{N} ( T(a) - H(at^{k}))P_{N} & \sim& 
G[a]^{N+l} E_{2}[a] \det P_{l}(T(a_0^{-1}) - H(a_0^{-1}t^{-1} )P_{l}\nn
\eq
as $N\to\iy$, where
\[
E_2[a] = \exp\Big( - \sum_{n =1}^{\infty} \log a_{2n}  -\frac{1}{2}\sum_{n=1}^\iy n  [\log a]_n^2 +\sum_{n=1}^\iy n [\log a]_{-n} [\log a]_n\Big).
\]
\item
Suppose that $k$ is a negative odd integer ($k=1-2l$, $l\ge 1$).
Then
\bq
\det P_{N} ( T(a) + H(at^{k}))P_{N}  &\sim&
  G[a]^{N+l} E_3[a] \det P_{l}(T(a^{-1}_0) + H(a^{-1}_0t ))P_{l}\nn
\eq
as $N\to\iy$, where
\[ 
E_3[a] = \exp\Big( -\log 2 + \sum_{n =1}^{\infty} \log a_{2n}  -\frac{1}{2}\sum_{n=1}^\iy n  [\log a]_n^2 +\sum_{n=1}^\iy n [\log a]_{-n} [\log a]_n\Big).
\]
\item 
We have
\bq
\det P_{N}  (T(a) + H(at^{k}))P_{N} &=&  0
\qquad \mbox{if } N\ge k\ge 2,\nn\\[1ex]
\det P_{N}  (T(a) - H(at^{k}))P_{N} &=&  0
\qquad \mbox{if } N\ge k\ge 1.\nn
\eq
\end{enumerate}
\end{theorem}
\begin{proof}
\underline{Case 1:}\
Consider the matrix $P_{N}  (T(a) \pm H(at^{k}))P_{N}.$ We observe that it is indeed the right bottom $N \times N$ corner of the $(N+l)\times(N+l)$ matrix
\[
A_{N}=P_{N+l} ( T(a) \pm  H(a))P_{N+l}.
\]
But this is the same the matrix $ Q_{l} A_{N} Q_{l}$. Using Jacobi's identity (\ref{f.J}) 
with $P=Q_l$, $Q=P_l$, we obtain
\[
\det(P_{N}  (T(a) \pm H(at^{k}))P_{N}) = \det(P_{l}A^{-1}_{N}P_{l})\cdot  (\det A_{N}) .
\]
Each of these last two factors can be computed asymptotically. 
For the second we use Theorem \ref{t4.3} with $M(a)=T(a)\pm H(a)$ and the results of Section \ref{s3} to conclude that  that $\det A_N$ is asymptotically $G[a]^{N+l} E_{1,\pm}[a]$. 

For the first we use the fact (\cite{BS}, Theorem 7.20) that the inverses of the finite sections $P_{N+l}  (T(a) \pm H(a))P_{N+l}$ converge strongly to the inverse of $M(a)=T(a) \pm H(a)$.
Notice that $M(a)$ is $T(a)$ plus a compact operator and that $M(a)$ is invertible.
The inverse equals $M(a)\iv=M(a_0\iv)T(a_-\iv)$. Hence  $\det\,P_{l}A^{-1}_NP_{l}$ 
converges to 
$$
\det P_{l} M(a_0\iv) T(a_-\iv) P_{l}=\det P_{l} M(a_0\iv) P_{l} T(a_-\iv) P_{l}=\det P_l M(a_0\iv) P_l.
$$
Here we use the basic fact that $T(a_-\iv)P_l=P_lT(a_-\iv)P_l$ and the normalization
$[a_-]_0=1$.

\underline{Case 2:}\
The proof is nearly the same only that we now use
$$
A_N=P_{N+l}(T(a)-H(a t\iv))P_{N+l}
$$
and $M(a)=T(a)-H(a t\iv)$.

\underline{Case 3:}\
Here an additional modification must be made. We consider
$$
A_N=P_{N+l}(T(a)+H(a t))P_{N+l}
$$
and using Jacobi's identity we can write
\[
\det(P_{N}  (T(a)-H(at^{k}))P_{N}) = \det(P_{l}A^{-1}_{N}P_{l})\cdot  (\det A_{N}) .
\]
Furthermore, we observe that
$$
M(a)=T(a)+H(a t) Q_1=(T(a)+H(at))R
$$
where 
$$
R=\mathrm{diag}(1/2,1,1,1,\dots).
$$
This leads to
$$
\det A_{N}=2\,\det (P_{N+l} M(a) P_{N+l}),
\quad
\det(P_{l}A^{-1}_{N}P_{l})=\frac{1}{2}\det(P_{l}(P_{N+l} M(a)P_{N+l})\iv P_{l}).
$$
The last determinant converges to
$$
\det(P_l M\iv(a)P_l)=\det(P_l M(a_0\iv) P_l)=\frac{1}{2}
\det(P_l(T(a_0\iv)+H(a_0\iv t))P_l).
$$
For this reason we get an additional factor $1/2$.

\underline{Case 4:}\
Observe that  $ P_{N} ( T(a) \pm  H(at^{k}))P_{N}$ is given by the matrix
\[ 
(a_{i-j} \pm a_{i+j-k+1})_{i,j=0\dots N-1}.
\] 
Thus the first column ($j = 0$) and the $k$th column ($j = k-1$) are given by
$a_{i}\pm a_{i-k+1}$ and $a_{i-k+1}\pm a_{i}$, respectively. Hence they are
either equal or the negative of each other. This settles the statements in the case
$k\ge 2$. The case $k=1$ with the ``minus'' is special. Then the first (= $k$th column)
equals zero.
 \end{proof}

\section{General $M(a)$}

In this section we consider the question of how general the operator $M(a)$ can be. Let us first make some preliminary basic observation. We will assume that $[M,\cS]$ is a compatible pair and that 
$\cS$ contains the trigonometric polynomials as a dense subset.

Let us write
\[ K(a) = M(a) - T(a). \]
Recall that the main property for compatible pairs implies the conditions (\ref{f.two}). The first one,
$M(ab) = M(a)M(b)$ for $b$ even, can be restated as 
\be\label{K.p}
 K(ab) = K(a)K(b) + T(a)K(b) + K(a)T(b) - H(a)H(b) 
\ee
whenever $b$ is even. The second one, $T(a)M(b)=M(ab)$ for $a \in \cS_{-}$, implies that 
for $a\in\cS_-$ we have $K(a)=0$  and $T(a)K(b) = K(ab)$ for any $b.$  

If we consider the matrix representation of $K(t^{n})$ with respect to the standard basis $e_{i} = \{\delta_{i,k}\}_{k=0}^\iy$ in $\ell^2$ ($i\ge 0$), then
\[ \langle K(t^{n})e_{i}, e_{j} \rangle =  \langle K(t^{n})e_{i}, T(t^{j})1 \rangle =  \langle T(t^{-j})K(t^{n})e_{i}, 1\rangle = \langle K(t^{n-j})e_{i}, 1\rangle = 0\]
when $j \geq n.$ In other words, $K(t^n)$ can have nonzero entries only in the first $n$ rows.
In particular, $K(t)$ is either zero or a rank one operator.

Moreover, since (for $n\ge 1$)
$$
K(t^{n+1})=K(t(t^{n}+t^{-n})),
$$
we can easily see from (\ref{K.p}) that once we have determined $K(t)$ we know $K(\cdot)$ for  arbitrary trigonometric polynomials. 

In what follows, let $e_0x^T$ with $x\in \ell^2$ stand for the rank one operator
$$
y\in \ell^2\mapsto e_0\langle y, x\rangle\in \ell^2
$$
We start with a proposition, which summarizes the main points we established so far.

\begin{proposition}
Let $[M,\cS]$ be a compatible pair, and assume that $\cS$ contains all trigonometric polynomials. Then
there exists $x\in \ell^2$ such that
\begin{enumerate}
\item[(i)] $M(t^{-n})=T(t^{-n})$ for $n\ge 0$,
\item[(ii)] $M((t+t\iv)^n)=(T(t+t\iv)+e_0x^T)^n$ for $n\ge 1.$
\end{enumerate}
\end{proposition}
\begin{proof}
The first assertion follows from the property $M(ab)=T(a)M(b)$ with $a=t^{-n}\in\cS_-$ and
$b=1$, noting that $M(1)=I$. The second property follows the multiplicative property in the case of even functions and from the fact that $M(t+t\iv)$ is of the form $T(t+t\iv)+e_0 x^T$. This last fact can be seen either from the remarks made above or from the equation
$$
T(t\iv)M(t+t\iv)=M(1+t^{-2})=I+T(t^{-2})=T(t\iv)T(t+t\iv),
$$
whence $T(t\iv) K(t+t\iv)=0$ and $K(t+t\iv)=e_0x^T$.
\end{proof}

Let $F\ell^1_1$ stand for the set of all functions $a$ defined on the unit circle such that its
Fourier coefficients satisfy
\[
\|a\|_{F\ell^1_1}:=\sum_{n\in\Z} (1+|n|) |a_n|<\iy.
\]
It is well known that $F\ell^1_1$ is a Banach algebra which is continuously embedded
into $C(\T)$ and has the set of trigonometric poplynomials as a dense subset.

\begin{theorem}
Let $\cS=F\ell^1_1$ and $x\in \ell^2$. Then the relations (i) and (ii) determine uniquely a well defined bounded linear operator for all trigonometric polynomials. Suppose that the operator norms of $K(t^{n})$ are uniformly bounded. Then in addition, the relations (i) and (ii) determine uniquely a well defined bounded linear operator for all $M:\cS\to \cL(\ell^2)$. Moreover, the pair $[M,\cS]$ is compatible.
\end{theorem}
\begin{proof}
It is obvious that (i) and (ii) determine a linear operator defined on the space $\cT$ of all  trigonometric
polynomials. We first want to show that
$$
M(abc)=T(a)M(b)M(c)
$$
holds for $a\in \cT_-$, $b\in\cT$, and $c\in \cT_0$. Since an arbitrary function in $\cT$ can be uniquely represented as a linear combination of $t^{-n}$ and $(t+t\iv)^{n}$, $n\ge 0$, it is not hard to see 
(see also (\ref{f.two})) that the only problem is to prove that $M(ab)=T(a)M(b)$ for $a\in \cT_-$, $b\in \cT_0$. We will consider
$a=t^{-n}$ and $b=t^m+t^{-m}$ and use induction on $n+m$. There is nothing to prove when $n=0$ or $m=0$. The case $n=m=1$ follows from
$$
T(t^{-1})M(t+t\iv)=T(t\iv)(T(t+t\iv)+e_0 x^T)=T(1+t^{-2})=M(1+t^{-2}).
$$
Now let $n,m\ge 1$ and $n+m>2$. Then
$$
T(t^{-n})M(t^m+t^{-m})=T(t^{-n+1})T(t^{-1})M(t^{m}+t^{-m})
$$
which by the induction arguments equals
\begin{align}
T(t^{-n+1})M(t^{m-1}+t^{-m-1}) &=T(t^{-n+1})M(t^{m-1}+t^{-m+1})+T(t^{-n+1})T(t^{ -m-1}-t^{-m+1})
\nn\\
&= M(t^{m-n}+t^{-m-n})+T(t^{-m-n}-t^{-n-m+2})
\nn\\
&=
M(t^{m-n}+t^{-n-m+2})\nn
\end{align}
as desired.

Finally, the remarks at the beginning of the section show that $K(t^{n}) = P_{n} K(t^{n})$ and since the trace norm of $P_{n}$ is $n$ we have that the trace norm of $K(t^{n})$ is bounded by some constant times $n$. Thus by our choice of  Banach algebra, our relations define a bounded linear operator for all function in $F\ell^1_1$.

\end{proof}

Our final remark is that it is easy to find examples that satisfy the hypothesis of our last theorem. In fact, we can take $x$ in our definition of $K$ to be the vector with a single one in the $j$th entry and zero otherwise.  In fact, our four concrete examples are when the first row of $K(t)$ is either the zero vector,
$\pm e_{0}$ or $e_{1}.$

%%%%%%%%%%%%%%%%%%%%%%%%%%%%%%%%%%%%%%%%%%%%%%%%%%%%%%%%%%%%%%%%%%%%%%%

\end{document}